\providecommand{\sectioninfo}{}
\documentclass{amsart}
\usepackage[utf8]{inputenc}
\usepackage{./preamble}

\title{Bicategorical Models of Classical Propositional Logic}
\author{Yuta Yamamoto}

\begin{document}

\maketitle

\begin{abstract}
F\"uhrmann and Pym constructed models of classical propositional logic in an order-enriched categorical setting, whose typical example is the category $\rel$ of sets and relations.
It is remarkable in that they are both non-degenerate and symmetric, i.e., free from the choices of the reduction strategy.

As a furter \textit{categorification} of this direction, we give bicategorical models of classical propositional logic that is also symmetric and non-degenerate.
Primal examples of our models include $\rel$, $\spancat$, and $\prof$, which shows that
we can construct models that are non-degenerate not only for $1$-cells but also for $2$-cells and the logical negations.

\end{abstract}

\section{Introduction}

While the categorical semantics of intuitionistic logic has been well studied, it is well known that a naive approach to classical logic causes difficulties.
This is related to the non-determinism of the process of cut elimination due to the symmetry of the rules of contraction and weakening.
Specifically, in the sequent calculus of classical propositional logic, it is known that for any two proofs of some formula, we can construct a proof that can be reduced to both of them (Lafont's example \cite{girard1989proofs}).
It is then challenging to construct non-degenerate semantics for symmetric contraction and weakening. While the well-known techniques related to the CPS transformation \cite{wadler03dual,selinger01control} provide non-degenerate models, they rely on a choice of particular asymmetric evaluation strategies.

In \cite{fuhrmann06order}, F\"uhremann and Pym constructed a semantics that would not depend on the arbitrary choice of such a particular strategy.
The key idea there was not to interpret reduction as strict equality, as is usually assumed, but as an order relation.
They interpret the entire system as an order-enriched category, including $\rel$ as a typical example.

However, the semantics of $\rel$ exhibits degeneracy in several other aspects.
First, while the structure of morphisms is non-degenerate, the structure of $2$-cells is degenerated.
Second, the interpretations of logical conjunction and disjunction coincide;
both of them are interpreted as the monoidal product in $\rel$, which is naturally regarded as a monoidal category with a dual.
Third, the logical negation is interpreted as the dual, which is collapsed to the identity functor.

To solve some of these degeneracies, let us consider the notion of bicategory, which is the categorification of the order-enrichment.
For instance, the bicategory of profunctors appears to admit non-degenerate interpretations for the first and third aspects.
The $2$-cells between profunctors are natural transformations.
Profunctors have a certain duality, that is, the opposite category.

Can we construct bicategorical semantics using the same approach as the symmetric order-enriched model for classical propositional logic?
In other words, can we devise a non-degenerate model that interprets propositions as objects, proofs as $1$-cells, and reductions as $2$-cells?
The challenge lies in how to handle reductions; it is necessary to assign $2$-cells as structures for each reduction.
As long as we treat reductions as orders, not only must we assign a $2$-cell for each reduction derivation, but also ensure that it is independent of the choice of derivation.
However, there is no reason for such coherence.
In fact, the naive attempt to obtain the $\prof$ model as the categorification of the $\rel$ model fails; we can easily find two derivations of some reduction which are interpreted as distinct $2$-cells.
To address this issue, it seems promising to consider a \textit{reduction-relevant} syntax that treats reductions not merely as orders but as structures, in addition to being proof-relevant.

In this study, we define the \textit{reduction-relevant} system {\ctwo} (a $2$-dimensional system of classical propositional logic) and construct symmetric and more non-degenerate bicategorical semantics.
Our models include $\rel$, $\spancat$, and $\prof$.
The main obstacle to the result is that, unlike usual type theories, syntactical substitution by variables does not correspond to the semantical composition of morphisms.
The reduction-relevant treatment of terms inevitably causes the notion of substitution splits into two: the usual substitution by variables and the \textit{cut}.
The syntax corresponding to the semantical composition is not the substitution but the cut.


\paragraph{Limitations and future work}
Although not pursued further in this paper, the bicategorical interpretation raises a natural question of whether this assignment of $2$-cells further satisfies certain ($3$-dimensional) equations.
Indeed, we can obtain {\cthree}, {\cfour}, and so on.
However, extending the system in this ad-hoc manner obscures the abstract structure behind them.
It is a natural question whether we can get something like {\cinfty} at once.

To achieve such high-dimensional extensions, considering type dependency (like \textit{homotopy type theory} \cite{hottbook}) seems to be a natural choice, but it faces some obstacles.
First, the system should be reduction-relevant, unlike the usual dependent type theory where the implicit computation of terms in a reduction-irrelevant way is assumed for type-checking.
Otherwise, it is likely to lead to degeneracy as \cite{herbelin05degeneracy} suggests.
Second, the corresponding semantics, i.e., the theory of ``$\infty$-profunctors'', has not been thoroughly examined.

Despite the challenges in this direction, there is a point worth mentioning. It is related to the interpretation of contradictions, which are interpreted in terms of hom-sets (more precisely, their generalization using presheavs and co-presheavs). This suggests one approach towards what is called \textit{directed type theory} \cite{licata11directed,riehl17type, cisinski23univalent}.

Another challenge to be addressed in the future is the semantical aspect of substitutions.
In the usual $\lambda$-calculus, given two terms, $M : A$ and $x : A \vdash N : B$, they can be ``composed'' in two ways:
$(\lambda x.N)M : B$ and $N[M/x]$.
We refer to the former as \textit{cut}, while the latter is called substitution.
In the cartesian closed category model, cuts and substitutions coincide.
On the other hand, in a symmetric system, like the one of ours, the semantical aspect of substitution is unclear.
In the symmetric systems, similar to the ordinary $\lambda$-calculus, one has the cut and the substitution for two (compatible) terms, and furthermore, one can consider the dual notion of substitution which we call co-substitution.
The two $\beta$ reductions in the symmetric systems can be understood as connecting cut and substitution (call-by-value) and cut and co-substitution (call-by-name), respectively.
While cut corresponds to the usual semantic composition, i.e., composition in the sense of profunctors, it seems that we have no straightforward semantics for substitution.


\section{Syntax of Classical Propositional Logic}

\sectioninfo

Since there are several possible variants in the syntax of {\ctwo}, we will provide two comments.
First, Our syntax is presented in a text style as in \cite{wadler03dual, curien00duality, lovas06classical} rather than in the graphical style as in \cite{fuhrmann06order}.
Second, Our syntax is presented in a polarity style; while in \cite{wadler03dual, curien00duality}, the ambient context is split into two parts, ours has a single context with the notion of a polarity which is the same as that of \cite{lovas06classical}.
    This is just a matter of convenience since the rule of \textit{exchange} is admissible.
We proceed by giving the raw syntax, forms of judgments, and rules.

\subsection{Raw syntax}
We have three kinds of notions of a ``term'' as in the literature, which we call term, co-term, and absurdity.
We fix a countable set of variables; we follow the convention and use $x$, $y$, {\ldots} for terms and $\alpha$, $\beta$, {\ldots} for co-terms.
We assume the following grammar:
\begin{alignat*}{2}
  A, B, C,  \ldots & ::=
    \tyTop \mid \tyBot \mid \tyAnd{A}{B} \mid \tyOr{A}{B} \mid \tyNot{A}
    &&\quad\text{(types)}\\
  L, M, N, \ldots & ::=
    x \mid \pbc{\alpha}{S} \mid
      \inTop \mid
      \inAnd{M}{N} \mid
      \elOrInl{M} \mid \elOrInr{M} \mid
      \elNot{K}
    &&\quad\text{(terms)}\\
  I, J, K, \ldots & ::=
    \alpha \mid \copbc{x}{S} \mid
      \inBot \mid
      \elAndFst{K} \mid \elAndSnd{K} \mid
      \inOr{J}{K} \mid
      \inNot{M}
    &&\quad\text{(co-terms)}\\
  S, T, U, \ldots & ::=
    \contra{M}{K}{A}
    &&\quad\text{(absurdity)}\\
  \pi, \sigma, \ldots & ::=
    \redRefl{M} \mid \redTrans{\pi}{\sigma} \mid
      \betamu{M}{x}{S}\mid
      \cobetamu{K}{x}{S}\mid \\
    &\quad
  \betaFst{M}{N}{K} \mid
  \betaSnd{M}{N}{K} \mid\\
    &\quad
  \betaInl{K}{L}{M} \mid
  \betaInr{K}{L}{N} \mid\\
    &\quad
  \betaNot{M}{K} \mid
  \congmu{\alpha}{\sigma} \mid
  \cocongmu{x}{\sigma} \mid
  \congContra{\pi}{\sigma} \mid\\
    &\quad
    \text{(congruences for logical connectives)} \cdots
    &&\quad\text{(reductions)}
%
\end{alignat*}
We omit congruence rules for term constructors for logical connectives, which can be easily verified to be added to the following discussion.
The dot notation, e.g., $\betamu{M}{x}{S}$, is to make it explicit that the right-hand side of the dot is \textit{bound} by the variable to the left of the dot.
It is as usual that substitutions should be done in a capture-avoiding way.


\subsection{Forms of judgment}
We assume the following six forms of judgment:
\begin{alignat*}{2}
  \termCtx{M}{+A}{\Gamma&} &&\quad\text{(terms)}\\
  \termCtx{K}{-A}{\Gamma&} &&\quad\text{(co-terms)}\\
  \termCtx{S}{\absurd}{\Gamma&} &&\quad\text{(absurdity)}\\
  \reductionCtx{\pi}{M}{N}{+A}{\Gamma&} &&\quad \text{(reductions)}\\
  \reductionCtx{\pi}{K}{L}{-A}{\Gamma&} &&\quad \text{(co-reductions)}\\
  \reductionCtx{\pi}{S}{T}{\absurd}{\Gamma&} &&\quad \text{(absurd reductions)}
\end{alignat*}


\subsection{Rules}
%

We list the rules in Figure 
\ref{fig:ctwo-rules-structural}, 
\ref{fig:ctwo-rules-structural-red}, 
\ref{fig:ctwo-rules-logical}, and
\ref{fig:ctwo-rules-logical-red} together with two comments.
First, we omit the ambient context $\Gamma$ except for (\ruleVR) and (\ruleVL).
For example, we write $\termCtx {M}{+B}{x: +A}$ to denote the term $\termCtx {M}{+B}{\Gamma, x: +A}$.
Second, in the rules (\ruleRefl) and (\ruleComp), we write, say, $\term{U}{\J}$ to denote an arbitrary term, co-term, or absurdity.

\begin{figure}[tb]
  \centering

  \AXC{$x: +A \in \Gamma$}
  \RL{(\ruleVR)}
  \UIC{$\term{x}{+A}$}
  \DP
  \quad %
  \AXC{$\alpha: -A\in\Gamma$}
  \RL{(\ruleVL)}
  \UIC{$\term{\alpha}{-A}$}
  \DP

  \vs
  \AXC{$\term{M}{+A}$}
  \AXC{$\term{K}{-A}$}
  \RL{(\ruleContra)}
  \BIC{$\term{\contra{M}{K}{A}}{\absurd}$}
  \DP

  \vs
  \AXC{$\termCtx{S}{\absurd}{\alpha: -A}$}
  \RL{(\ruleRI)}
  \UIC{$\term{\pbc{\alpha}{S}}{+A}$}
  \DP
  \quad %
  \vs
  \AXC{$\termCtx{S}{\absurd}{x: +A}$}
  \RL{(\ruleLI)}
  \UIC{$\term{\copbc{x}{S}}{-A}$}
  \DP

  \caption{Structural terms}
  \label{fig:ctwo-rules-structural}
\end{figure}
\begin{figure}[tb]
  \centering

  \AXC{$\term{K}{-A}$}
  \AXC{$\termCtx{S}{\absurd}{\alpha: -A}$}
  \RL{(\ruleBetaRI)}
  \BIC{$\reduction{\betamu{K}{\alpha}{S}} { \contra{\pbc{\alpha}{S}}{K}{A}} { \subst{K}{\alpha}{S}} {\absurd}$}
  \DP

  \vs
  \AXC{$\term{M}{+A}$}
  \AXC{$\termCtx{S}{\absurd}{x: +A}$}
  \RL{(\ruleBetaLI)}
  \BIC{$\reduction{\cobetamu{M}{x}{S}} { \contra{M}{\copbc{x}{S}}{A}} { \subst{M}{x}{S}} {\absurd}$}
  \DP

  \vs
  \AXC{$\term{U}{\J}$}
  \RL{(\ruleRefl)}
  \UIC{$\reduction{\redRefl{U}}{U}{U}{\J}$}
  \DP

  \vs
  \AXC{$\reduction{\pi}{U}{V}{\J}$}
  \AXC{$\reduction{\sigma}{V}{W}{\J}$}
  \RL{(\ruleComp)}
  \BIC{$\reduction{
    \redTrans{\pi}{\sigma}}{U}{W}{\J}$}
  \DP


  \vs
  \AXC{$\reductionCtx{\sigma}{S}{T}{\absurd}{\alpha: -A}$}
  \RL{(\ruleCongRI)}
  \UIC{
    $\reduction{\congmu{\alpha}{\sigma}}
    { \pbc{\alpha}{S}}{ \pbc{\alpha}{T}}
    {+A}$}
  \DP

  \vs
  \AXC{$\reductionCtx{\sigma}{S}{T}{\absurd}{x: +A}$}
  \RL{(\ruleCongLI)}
  \UIC{
    $\reduction{\cocongmu{x}{\sigma}}
    { \copbc{x}{S}}{ \copbc{x}{T}}
    {-A}$}
  \DP

  \vs
  \AXC{$\reduction{\pi}{M}{N}{+A}$}
  \AXC{$\reduction{\sigma}{K}{L}{-A}$}
  \RL{(\ruleCongContra)}
  \BIC{$\reduction{\congContra{\pi}{\sigma}}
  {\contra{M}{K}{A}}
  {\contra{N}{L}{A}}
  {\absurd}$}
  \DP

%

  \caption{Reductions for structural terms}
  \label{fig:ctwo-rules-structural-red}
\end{figure}

\begin{figure}[tb]
  \begin{center}
  \fbox{$\tyTop$}
  \AXC{}
  \RL{($+\tyTop$)}
  \UIC{$\term{\inTop}{+\tyTop}$}
  \DP
  \quad
  \fbox{$\tyBot$}
  \AXC{}
  \RL{($-\tyBot$)}
  \UIC{$\term{\inBot}{-\tyBot}$}
  \DP

  \vs
    \fbox{$\tyAnd{}{}$}
  \AXC{$\tmPos{M}{A}$}
  \AXC{$\tmPos{N}{B}$}
  \RL{($+\tyAnd{}{}$)}
  \BIC{$\tmPos{\inAnd{M}{N}}{\tyAnd{A}{B}}$}
  \DP

  \vs
    \AXC{$\term{K}{-A}$}
    \RL{($-\land_1$)}
    \UIC{$\term{\elAndFst{K}}{-\tyAnd{A}{B}}$}
  \DP
    \AXC{$\term{K}{-B}$}
    \RL{($-\land_2$)}
    \UIC{$\term{\elAndSnd{K}}{-\tyAnd{A}{B}}$}
  \DP


  \vs
    \fbox{$\tyOr{}{}$}
  \AXC{$\term{K}{-A}$}
  \AXC{$\term{L}{-B}$}
  \RL{($-\tyOr{}{}$)}
  \BIC{$\term{\inOr{K}{L}}{-\tyOr{A}{B}}$}
  \DP

  \vs
    \AXC{$\term{M}{+A}$}
    \RL{($+\tyOr{}{}_1$)}
    \UIC{$\term{\elOrInl{M}}{+\tyOr{A}{B}}$}
  \DP
    \AXC{$\term{M}{-B}$}
    \RL{($+\tyOr{}{}_2$)}
    \UIC{$\term{\elOrInr{M}}{-\tyOr{A}{B}}$}
  \DP

  \vs
    \fbox{$\tyNot{}$}
  \AXC{$\term{K}{-A}$}
  \RL{($+\tyNot{}$)}
  \UIC{$\term{\inNot{K}}{+\tyNot{A}}$}
  \DP
\quad
  \AXC{$\term{M}{+A}$}
  \RL{($-\tyNot{}$)}
  \UIC{$\term{\elNot{M}}{-\tyNot{A}}$}
  \DP
  \end{center}

  \caption{Terms for logical connectives}
  \label{fig:ctwo-rules-logical}
\end{figure}
\begin{figure}[tb]
  \begin{center}
  \AXC{$\tmPos{M}{A}$}
  \AXC{$\tmPos{N}{B}$}
  \AXC{$\tmNeg{K}{A}$}
  \RL{(\ruleBetaFst)}
  \TIC{$\reduction{\betaFst{M}{N}{K}}
  {\contra{\inAnd{M}{N}}{\elAndFst{K}}{\tyAnd{A}{B}} }
  {\contra{M}{K}{A}}
    {\absurd}$}
  \DP

  \vs
  \AXC{$\term{M}{+A}$}
  \AXC{$\term{N}{+B}$}
  \AXC{$\term{K}{-B}$}
  \RL{(\ruleBetaSnd)}
  \TIC{$\reduction{\betaSnd{M}{N}{K}}
  {\contra{\inAnd{M}{N}}{\elAndSnd{K}}{\tyAnd{A}{B}} }
  {\contra{N}{K}{B}}
    {\absurd}$}
  \DP

  \vs
  \AXC{$\term{K}{-A}$}
  \AXC{$\term{L}{-B}$}
  \AXC{$\term{M}{+A}$}
  \RL{(\ruleBetaInl)}
  \TIC{$\reduction{\betaInl{K}{L}{M}}
  {\contra{\elOrInl{M}}{\inOr{K}{L}}{\tyOr{A}{B}} }
  {\contra{K}{M}{A}}
    {\absurd}$}
  \DP

  \vs
  \AXC{$\term{K}{-A}$}
  \AXC{$\term{L}{-B}$}
  \AXC{$\term{M}{+B}$}
  \RL{(\ruleBetaInr)}
  \TIC{$\reduction{\betaInr{K}{L}{M}}
  {\contra{\elOrInr{M}}{\inOr{K}{L}}{\tyOr{A}{B}} }
  {\contra{L}{M}{B}}
    {\absurd}$}
  \DP

  \vs
  \AXC{$\term{M}{+A}$}
  \AXC{$\term{K}{-A}$}
  \RL{(\ruleBetaNeg)}
  \BIC{$\reduction{\betaNot{M}{K}}
  {\contra{\elNot{M}}{\inNot{K}}{\tyNot{A}} }
  {\contra{M}{K}{A}}
    {\absurd}$}
  \DP

  \end{center}

  \caption{Reductions for logical connectives}
  \label{fig:ctwo-rules-logical-red}
\end{figure}


\section{Semantics}

\sectioninfo

So far we have defined our syntax. Now we construct the bicategorical model.
More concretely, we work on the bicategory $\prof$.
It is in some sense quite natural in that the following construction can also be applied to $\rel$ and $\spancat$ by reading relations and spans as certain profunctors.

\subsection{The bicategory of profunctors}
For completeness and in order to fix notations, we begin by recalling the formal definitions of profunctor (also called distributor).
The canonical reference is \cite{benabou2000distributors}.
\begin{definition}
  Let $\A$ and $\B$ be (small) categories.
  A profunctor $P$ from $\A$ to $\B$, written as $P: \A \proTo \B$,
  is a functor $P: \A \times \op{\B} \to\set$.
  It is equivalently a functor from $\A$ to the presheaf category on $\B$.
\end{definition}

For convenience, we define its multivariable variant.
It is similar to functor in that it can consider many-inputs, but profunctor can also take ``many-outputs'':
\begin{definition}
  Let $n,m\geq 0$ be natural numbers and $\A_1$,...,$\A_n$, $\B_1$,...,$\B_m$ be categories.
  A (multi-variable) profunctor
  $P: \A_1,...,\A_n \proTo \B_1,...,\B_m$
  is a profunctor
  $P: \A_1\times...\times\A_n \proTo \B_1\times...\times\B_m$.
  When $n=m=1$, it is a profunctor in a usual sense.
  When $n=0$ and $m=1$, it is a presheaf on $\B_1$.
  When $n=1$ and $m=0$, it is a co-presheaf on $\A_1$.
  When $n=m=0$, it is just a set.
\end{definition}

\begin{notation}
First, we use letters $\Gamma$ and $\Delta$ to denote finite lists of categories.
We will often be sloppy about the distinction between the list itself and the category defined by its direct product.
For example, we write $P: \Gamma \times \op{\Delta} \to \set$ to denote a functor which is equivalent to a multivariable profunctor $P:\Gamma \proTo \Delta$.
Second, we write $P, Q:\Gamma \proTo \Delta$ to denote two profunctors $P:\Gamma \proTo \Delta$ and $Q:\Gamma \proTo \Delta$.
\end{notation}

\begin{definition}
  Let $\A$ be a category. An identity profunctor $\id{\A}:\A\proTo\A$ is defined to be a functor
  $\hom_\A: \A\times\op\A\to\set$.
\end{definition}

\begin{definition}
  Let $\A$, $\B$, and $\C$ be categories and $P: \A \proTo \B$ and $Q: \B \proTo \C$ be profunctors.
  Their composite profunctor $Q\circ P:\A \proTo \C$ is defined to be the functor
  \begin{alignat*}{2}
    x,z &\mapsto \coend{y\in \B}{P(x, y)\times Q(y, z)}:
    \A \times \op\C &\to\set
  \end{alignat*}
\end{definition}
It is useful to have a multivariable variant of composition of profunctors:
\begin{definition}
  Let $\Gamma$ and $\Delta$ be finite lists of categories.
  Let $P:\Gamma\proTo\Delta,\A$ and $Q:\A,\Gamma\proTo\Delta$ be profunctors.
  Their composite profunctor $Q\circ_\A P:\Gamma\proTo\Delta$ is defined to be the functor
  \begin{alignat*}{2}
    \gamma,\delta &\mapsto \coend{x\in \A}{P(\gamma,\delta, x)\times Q(x, \gamma, \delta)}:
    \Gamma\times \op\Delta &\to\set
  \end{alignat*}
\end{definition}
It can also be defined by Kan extensions as in \cite{benabou2000distributors}, but the definition using coends is more suitable for the following argument because of its symmetry.

\begin{definition}
  Let $P,Q:\Gamma\proTo\Delta$ be parallel profunctors.
  A $2$-cell between them is a natural transformation between them regarded as $\set$-valued functors.
\end{definition}
These define a bicategory $\prof$ whose objects are (small) categories and $1$-cells are profunctors.

\subsection{The interpretation}
We use $\eval{-}$ to denote the interpretation function but we often blur the distinction between syntax and semantics in order not to clutter things up with brackets.
We separate the positive and negative parts of the context, which we denote by $+$ and $-$ subscripts, respectively.
We intend to interpret the syntax by the following schema:
\begin{itemize}
  \item type $A$ as a category;
  \item term $\termCtx{M}{+A}{\Gamma}$ as a profunctor
    $\Gamma_+ \proTo \Gamma_-, A$;
  \item coterm $\termCtx{K}{-A}{\Gamma}$ as a profunctor
    $\Gamma_+, A \proTo \Gamma_-$;
  \item absurdity $\termCtx{S}{\absurd}{\Gamma}$ as a profunctor
    $\Gamma_+ \proTo \Gamma_-$; and
  \item reduction $\reductionCtx{\pi}{V}{W}{\J}{\Gamma}$ as a natural transformation $\pi: V\To W$.
\end{itemize}
The main difficulty lies in the construction of structural $\beta$, i.e., (\ruleBetaRI) and (\ruleBetaLI).
We proceed by giving the interpretations for types, structural terms, structural reductions except for $\beta$, logical terms and reductions, and structural $\beta$.
\paragraph{Types}
\begin{alignat*}{2}
  \tyTop &:= 1\\
  \tyBot &:= 1\\
  \tyAnd{A}{B} &:= A\times B\\
  \tyOr{A}{B} &:= A\times B\\
  \tyNot{A} &:= \op A
\end{alignat*}

%
%

\paragraph{Structural terms and reductions except for $\beta$}
\begin{itemize}
  \item Variable $\termCtx{x}{\A}{\Gamma}$ is mapped to the Yoneda embedding (pre-composed with the projection)
    $$\Gamma\to \presh{\A}$$

  \item Cut $\termCtx {\contra {M}{K}{A}}{\absurd}{\Gamma}$ is interpreted by the composition of profunctors, i.e.,
    $K \circ_A M: \Gamma_+ \proTo \Gamma_-$.

  \item The above construction is natural in $M$ and $K$, which defines the interpretation of $\contra{\sigma}{\pi}{}$.

  \item The interpretations of $\pbc{\alpha}{S}$, $\copbc{x}{S}$, $\pbc{\alpha}{\sigma}$ and $\copbc{x}{\sigma}$ are trivial.

  \item Reflexivity and transitivity are interpreted as identities and compositions of $2$-cells.
\end{itemize}
\paragraph{Logical terms and reductions} 
We first work on the connective $\land$:
\begin{itemize}
  \item
    Given presheaves $\term{M}{+A}$ and $\term{N}{+B}$,
    a presheaf $\term{\inAnd {M}{N}}{+\tyAnd{A}{B}}$
    is defined by $\inAnd{M}{N}(\pair{x, y}) := M(x)\times N(y)$.
  \item
    Given a copresheaf $\term{K}{-A}$,
    a copresheaf $\term{\elAndFst{K}}{-\tyAnd{A}{B}}$
    is defined by $\elAndFst{K}(\pair{x, y}) := K(x)$.
  \item
    $\reduction{\betaFst{M}{N}{K}}
      {\contra{\inAnd{M}{N}}{\elAndFst{K}}{\tyAnd{A}{B}} }
      {\contra{M}{K}{A}}
      {\absurd}$ is defined to be the canonical $2$-cell
    \begin{alignat*}{2}
    {\contra{\inAnd{M}{N}}{\elAndFst{K}}{\tyAnd{A}{B}} }
      &= \coend{x, y}{\inAnd{M}{N}(x, y)\times \elAndFst{K}(x, y)}\\
      &= \coend{x, y}{(M(x)\times N(y))\times K(x)}\\
      &\to \coend{x}{M(x)\times K(x)}\\
      &= {\contra{M}{K}{A}}.
    \end{alignat*}


\end{itemize}
Similarly, the connective $\neg$ can be constructed by the fact that a presheaf on a category is equivalent to a copresheaf on its opposite category.
Other cases can be done in a similar manner.
\paragraph{Structural $\beta$}
By symmetry, we focus on (\ruleBetaRI), i.e.,
$\reduction{\betamu{K}{\alpha}{S}}
  {\contra{\pbc{\alpha}{S}}{K}{A}}
  {\subst{K}{\alpha}{S}}
  {\absurd}$.
We proceed by induction on the structure of $S$.
Since $S$ is of the form $\contra{N}{L}{}$, the induction will be done also on the structure of $N$ and $L$, which leads to the idea of the following \textit{delayed substitution} and generalized $\beta$:

\begin{definition}
  Let $K:-A$ and $\termCtx{V}{\J}{\alpha:-A}$ be terms.
  A delayed substitution $\substex{K}{\alpha}{V}:\J$ is defined as:
  \begin{alignat*}{2}
    \substex{K}{\alpha}{S} &:= \contra{\pbc{\alpha}{S}}{K}{A} :\absurd\\
    \substex{K}{\alpha}{N} &:= \pbc{\beta}{(\substex{K}{\alpha}{\contra{N}{\beta}{B}})} :+B\\
    \substex{K}{\alpha}{L} &:= \pbc{y}{(\substex{K}{\alpha}{\contra{y}{L}{B}})} :-B
  \end{alignat*}
  Dually, $\substex{M}{x}{V} :\J$ is defined.
\end{definition}

Now we generalize the structural $\beta$ as:
\begin{alignat*}{2}
  \reduction{\betamu{K}{\alpha}{V}}
    {\substex{K}{\alpha}{S}}
    {\subst{K}{\alpha}{S}}
    {\J}
\end{alignat*}
It is clear that the original (\ruleBetaRI) is its special case.
We will construct the generalized $\beta$ and obtain the original one as a corollary.
\begin{construction}[Construction of the interpretation for generalized structural $\beta$]
We proceed by structural induction on $V$:
\begin{itemize}
  \newcommand{\localTermV}[0]{\contra{L}{K}{B}}
  \item case where $V=\localTermV$:\\
    $\redTm
      {\betamup {M}{x}{\localTermV}}
      {\substex{M}{x}{\localTermV}}
      {\subst{M}{x}{\localTermV}}$
    is the composite
    \begin{alignat*}{2}
      \eval{\substex{M}{x}{\contra{L}{K}{B}} }
        &= \coend{x\in A}{M(x)\times \eval{\contra{L}{K}{B}}(x)}\\
        &= \coend{x\in A}{M(x)\times
          \bigparen{\coend{y\in B}{L(x, y)\times K(x, y)}} }\\
        &\to \coend{y\in B}{
          \bigparen{\coend{x\in A}{M(x)\times L(x, y)}}\times
          \bigparen{\coend{x\in A}{M(x)\times K(x, y)}} }\\
        &= \coend{y\in B}{
          \eval{\substex{M}{x}{L}}(y)\times
          \eval{\substex{M}{x}{K}}(y)}\\
        &\to \coend{y\in B}{
          \eval{\subst{M}{x}{L}}(y)\times
          \eval{\subst{M}{x}{K}}(y)}\\
        &= \eval{\contra{\subst{M}{x}{L}}{\subst{M}{x}{K}}{B}}\\
        &= \eval{\subst{M}{x}{\contra{L}{K}{B}} }
    \end{alignat*}
  \renewcommand{\localTermV}[0]{x}
  \item case where $V=\localTermV:+A$:\\
    $\reduction
      {\betamup {M}{x}{\localTermV}}
      {\substex{M}{x}{\localTermV}}
      {\subst{M}{x}{\localTermV}}
      {+A}$
    is the composite (which is actually an isomorphism by the ninja Yoneda)
    \begin{alignat*}{2}
      \eval{\substex{M}{x}{\localTermV}} 
        &= \eval{ \pbc{x'}{ \contra{M}{\pbc{x}{\contra{x}{x'}{A}} }{A}} }\\
        &= x'\mapsto \coend{x}{M(x)\times A(x', x)}\\
        &\to M\\
        &= \eval{\subst{M}{x}{\localTermV}} 
    \end{alignat*}
  \renewcommand{\localTermV}[0]{y}
  \item case where $V=\localTermV:+B (\neq x)$:\\
    $\reduction
      {\betamup {M}{x}{\localTermV}}
      {\substex{M}{x}{\localTermV}}
      {\subst{M}{x}{\localTermV}}
      {+B}$
    is the composite
    \begin{alignat*}{2}
      \eval{\substex{M}{x}{\localTermV}} 
        &= \eval{ \pbc{x'}{ \contra{M}{\pbc{x}{\contra{y}{x'}{A}} }{A}} }\\
        &= x'\mapsto \coend{x}{M(x)\times A(x', y)}\\
        &\to \y y\\
        &= \eval{y}\\
        &= \eval{\subst{M}{x}{\localTermV}} 
    \end{alignat*}
  \renewcommand{\localTermV}[0]{\pbc{\alpha}{S}}
  \item case where $V=\localTermV :+B$:\\
    $\reduction
      {\betamup {M}{x}{\localTermV}}
      {\substex{M}{x}{\localTermV}}
      {\subst{M}{x}{\localTermV}}
      {+B}$
    is the composite
    \begin{alignat*}{2}
      \eval{\substex{M}{x}{\localTermV}} 
        &= \alpha\mapsto \eval{ \substex{M}{x}{S}}\\
        &\to \alpha\mapsto \eval{ \subst{M}{x}{S}}\\
        &= \eval{\subst{M}{x}{\localTermV}} 
    \end{alignat*}
  \renewcommand{\localTermV}[0]{\inAnd{N_1}{N_2}}
  \item case where $V$ is an introduction term, e.g.,
    $V=\localTermV :+B_1 \times B_2$:\\
    $\reduction
      {\betamup {M}{x}{\localTermV}}
      {\substex{M}{x}{\localTermV}}
      {\subst{M}{x}{\localTermV}}
      {+B_1 \times B_2}$
    is the composite
    \begin{alignat*}{2}
      \eval{\substex{M}{x}{\localTermV}} 
        &= z\mapsto \coend{x}{M(x)\times \inAnd{N_1}{N_2}(x, z)}\\
        &= \pair{y_1, y_2}\mapsto \coend{x}{M(x)\times N_1(x, y_1)\times N_2(x, y_2)}\\
        &\to \pair{y_1, y_2}\mapsto
          \bigparen{\coend{x}{M(x)\times N_1(x, y_1)}} \times
          \bigparen{\coend{x}{M(x)\times N_2(x, y_2)}}\\
        &= \pair{y_1, y_2}\mapsto
          \substex{M}{x}{N_1}(y_1)\times
          \substex{M}{x}{N_2}(y_2)\\
        &\to \pair{y_1, y_2}\mapsto
          \subst{M}{x}{N_1}(y_1)\times
          \subst{M}{x}{N_2}(y_2)\\
        &= \eval{\subst{M}{x}{\localTermV}} 
    \end{alignat*}
  \renewcommand{\localTermV}[0]{\elAndFst{K}}
  \item case where $V$ is an elimination term, e.g., $V=\localTermV :-B_1 \times B_2$:\\
    $\reduction
      {\betamup {M}{x}{\localTermV}}
      {\substex{M}{x}{\localTermV}}
      {\subst{M}{x}{\localTermV}}
      {-B_1 \times B_2}$
    is the composite
    \begin{alignat*}{2}
      \eval{\substex{M}{x}{\localTermV}} 
        &\iso z\mapsto \coend{x}{M(x)\times (\elAndFst{K}(z))}\\
        &= \pair{y_1, y_2}\mapsto \coend{x}{M(x)\times K(y_1)}\\
        &= \pair{y_1, y_2}\mapsto \substex{M}{x}{K}(y_1)\\
        &\to \pair{y_1, y_2}\mapsto \subst{M}{x}{K}(y_1)\\
        &= \eval{\subst{M}{x}{\localTermV}} 
    \end{alignat*}
\end{itemize}
\end{construction}
Since this generalized one can be defined via the original one,
  e.g., $\betamu{K}{\alpha}{M} := \congmu{\beta}{\betamu{K}{\alpha}{\contra{M}{\beta}{}} }$,
  one might wonder whether their interpretations coincide.
  They indeed coincide and this kind of higher dimensional equalities should be studied in the future.


\bibliography{./bicat-model}

\end{document}